\documentclass[10pt, reqno]{amsart}
\usepackage{amsmath}
\usepackage{amssymb}
\usepackage{amsaddr}
\usepackage{enumerate}
\usepackage{amsbsy}
\usepackage{amsfonts}
\usepackage{color}
\usepackage{upgreek}

\newtheorem{thm}{Theorem}[section]
\newtheorem{lem}[thm]{Lemma}

\theoremstyle{remark}
\newtheorem{rem}[thm]{Remark}

\newcommand{\les}{\lesssim}

	\title[Charge conjugation approach]{Charge conjugation approach to scattering for the Hartree type Dirac equations with chirality}

\author[Y. Cho]{Yonggeun Cho}
\address{Department of Mathematics, and Institute of Pure and Applied Mathematics, Jeonbuk National University, Jeonju 54896, Republic of Korea}
 \email{changocho@jbnu.ac.kr}

\author[S. Hong]{Seokchang Hong}
    \address{Department of Mathematics, Chung-Ang University, Seoul 06974, Republic of Korea}
    \email{seokchangh11@cau.ac.kr}

	\author[T. Ozawa]{Tohru Ozawa}
	\address{Department of Applied Physics, Waseda University, 3-4-1, Okubo, Shinjuku-ku, Tokyo, 169-8555, Japan}
	\email{txozawa@waseda.jp}

\begin{document}

	\thanks{2020 {\it Mathematics Subject Classification.} M35Q55, 35Q40.}
	\thanks{{\it Key words and phrases.} Hartree-type Dirac equation, global well-posedness, scattering, charge conjugation, chirality}

\begin{abstract}
	We study the Cauchy problems for the Hartree-type nonlinear Dirac equations with Yukawa-type potential derived from pseudoscalar field. We establish scattering for large data but with a relatively small part of initial data associated with charge conjugation by exploiting null structure induced by chiral operator.
	\end{abstract}

		\maketitle

\section{Introduction}
The purpose of this paper is to investigate scattering property for the cubic Dirac equation with the Hartree-type nonlinearity in $\mathbb R^{1+3}$ given by
\begin{align}\label{main-eq}
\left\{
\begin{array}{l}
	-i\gamma^\mu\partial_\mu\psi + m\psi = [V_b*(\overline\psi \gamma^5\psi)]\gamma^5\psi \;\;\mbox{in}\;\;\mathbb R^{1 + 3}, \\
	\psi|_{t=0}:= \psi_0.
\end{array}\right.	
\end{align}
Here $m  > 0$, $*$ is the convolution in $\mathbb R^3$, $\overline{\psi} = \psi^\dagger \gamma^0$, $\psi^\dagger$ is the transpose of complex conjugate of $\psi$, and  $\gamma^5 = i\gamma^0\gamma^1\gamma^2\gamma^3$. The potential $V_b$ is the Yukawa potential interacting between nucleon and meson and is given by
$$
V_b(x) := \frac{g_0^2}{4\pi} \frac{e^{-b|x|}}{|x|},
$$
where $g_0 \in \mathbb R$ and $b > 0$ are given physical parameters.
The gamma matrices $\gamma^\mu\in\mathbb C^{4\times4}\;(\mu = 0, 1, 2, 3)$ are given by
\begin{align*}
\gamma^0 = \begin{bmatrix}
 	I_{2\times2} & \mathbf0 \\ \mathbf0 & -I_{2\times2}
 \end{bmatrix}, \ \gamma^j = \begin{bmatrix}
 	\mathbf 0 & \sigma^j \\ -\sigma^j & \mathbf0
 \end{bmatrix}	
\end{align*}
with the Pauli matrices $\sigma^j\in\mathbb C^{2\times2} \;(j = 1, 2, 3)$ given by
\begin{align*}
\sigma^1 = \begin{bmatrix}
 	0 & 1 \\ 1 & 0
 \end{bmatrix}, \ \sigma^2 = \begin{bmatrix}
 	0 & -i \\ i & 0
 \end{bmatrix}, \  \sigma^3=\begin{bmatrix}
 	1 & 0 \\ 0 & -1
 \end{bmatrix}.	
\end{align*}

In this paper we establish global well-posedness and scattering of solutions to the system \eqref{main-eq} for small charge conjugation in the scaling critical space which has extra weighted regularity in the angular variables. To be more precise, we let $\Omega_{ij}=x_i\partial_j-x_j\partial_i$ be the infinitesimal generators of the rotations on $\mathbb R^3$ and let $\Delta_{\mathbb S^2} = \sum_{1 \le i < j \le 3}\Omega_{ij}^2$ be the Laplace-Beltrami operator on the unit sphere $\mathbb S^2\subset\mathbb R^3$. Then we can define the fractional power of angular derivative by $\Lambda_{\mathbb S^2}^\sigma = (1-\Delta_{\mathbb S^2})^\frac\sigma2$, which will be treated concretely below, and define angularly regular space $L_x^{2, \sigma}$ space by $\Lambda_{\mathbb S^2}^{-\sigma}L_x^2$ and its norm by $\|f\|_{L_x^{2, \sigma}} := \|\Lambda_{\mathbb S^2}^\sigma f\|_{L_x^2}$.
Now we state the main theorem:

\begin{thm}\label{gwp}
Let $\sigma > 0$ and $\theta \in \{+, -\}$. Then there exists $\mathtt a = \mathtt a(\|\psi_0\|_{L^{2, \sigma}_x(\mathbb R^3)}) > 0$ such that for all initial data $\psi_0 \in L_x^{2, \sigma}$ satisfying
\begin{align}\label{majorana}
\|\psi_0 + i\gamma^2\psi_0^*\|_{L^{2, \sigma}_x(\mathbb R^3)} \le \mathtt a,
\end{align}
the Cauchy problem \eqref{main-eq} is globally well-posed and solution $\psi$ scatters in $L_x^{2, \sigma}$ to free solutions as $t\rightarrow\pm\infty$.
\end{thm}
\noindent Here $\psi_0^*$ denotes the complex conjugate of $\psi_0$. The transformed spinor $i\gamma^2\psi_0^*$ is referred to as a charge conjugation of $\psi_0$, which will be discussed below in detail.
The condition \eqref{majorana} was observed in \cite{candyherr1, chohlee1} for the equations with identity matrix in place of $\gamma^5$. It is regarded as a perturbation of Majorana condition studied in \cite{majorana, chagla, ozya}.

One may observe that the equation \eqref{main-eq} can be derived by decoupling the following Dirac-Klein-Gordon system
\begin{align}\label{dkg}
\begin{aligned}
&(-i\gamma^\mu\partial_\mu + m)\psi = g_0\phi i\gamma^5\psi, \\
&(\partial_t^2-\Delta+M^2)\phi = -g_0\overline{\psi}i\gamma^5\psi.
\end{aligned}
\end{align}
The system \eqref{dkg} conserves parity and is Lorentz covariant. The matrix $\gamma^5$ was chosen for the right-hand side to be a pseudoscalar.  For details, see Ch. 10 of \cite{bjor}.

Let us assume that the pseudoscalar field $\phi$ is a standing wave, i.e., $\phi(t,x)=e^{i\lambda t}f(x)$ with $M > |\lambda|$. Then the Klein-Gordon part of \eqref{dkg} becomes
\begin{align}\label{kg}
(- \Delta  + M^2 -\lambda^2)\phi = -g_0\overline{\psi}i\gamma^5\psi.
\end{align}
Then we put \eqref{kg} into the Dirac part of \eqref{dkg} and then a spinor field $\psi$ gives the desired equation.

If $\gamma^5$ is replaced by the identity (in this case $\phi$ is a scalar field), then \eqref{main-eq} and \eqref{dkg} have been extensively studied \cite{chagla, danfoselb, behe, candyherr, candyherr1, tes, tes1, cyang, geosha, cholee, chohlee, choozlee}. Our work is motivated from these works and is concerned with the relation between charge conjugation and chirality (see \eqref{charge-c-chirality} below).

\begin{rem}[Chirality] The gamma matrix $\gamma^5$ represents the chirality of spinors.
Let $\psi_R = \frac12(1 + \gamma^5)\psi$ and $\psi_L = \frac12(1-\gamma^5)\psi$. Then $\psi = \psi_R + \psi_L$. Using $\gamma^5\gamma^5 = 1$, $\gamma^5\psi_R = \psi_R$ and $\gamma^5\psi_L = -\psi_L$, which lead us to the chirality. Hence $(\psi_R)_R = \psi_R, (\psi_L)_L = \psi_L$ and $(\psi_R)_L = (\psi_L)_R = 0$. Since $\gamma^5\gamma^\mu \gamma^\mu\gamma^5 = 0$ for $\mu = 0, 1, 2, 3$, $\overline{\psi_R}\psi_R = \overline{\psi_L}\psi_L = 0$ and
$\gamma^5(-i\gamma^\mu\partial_\mu)\psi) = i\gamma^\mu\partial_\mu \gamma^5\psi$.
Hence we conclude that the solution $\psi$ to \eqref{main-eq} satisfies
\begin{align}\begin{aligned}\label{chiral-form}
-i\gamma^\mu\partial_\mu \psi_R  &= -M\psi_L + [V_b * (\overline\psi_L\psi_R - \overline\psi_R\psi_L)]\psi_L,\\
-i\gamma^\mu\partial_\mu \psi_L  &= -M\psi_R + [V_b * (\overline\psi_L\psi_R - \overline\psi_R\psi_L)]\psi_R,\\
\psi_R(0) &= \psi_{0, R},\;\;\psi_R(0) = \psi_{0, L}.
\end{aligned}\end{align}
The new system \eqref{chiral-form} shows that the right$($left$)$-handed chiral field $\psi_R(\psi_L)$ can evolve even though starting off as a completely left$($right$)$-handed chiral state due to the influence of left$($right$)$-handed chiral field. This shows that the chirality is not conserved and a smallness on both data $\psi_{0, R}$ and $\psi_{0, L}$ rather than a partial smallness is necessary for the global well-posedness unlike the charge conjugation as stated in Theorem \ref{gwp}.
\end{rem}
\begin{rem}
If $b = 0$, then \eqref{main-eq} becomes scattering-critical equation due to the $\log$-blowup nature. Hence a non-scattering is plausible in $L_x^2$ space unless the Majorana condition appears. In fact, one may consider scattering states satisfying coercivity condition guaranteeing the full time decay $t^{-\frac32}$ (see \cite{choozxia}) as in \cite{choozlee, chohlee}, to which any solution of \eqref{main-eq} dose not converge in $L_x^2$. Also one may try to remove the $\log$-blowup by modifying phases. It would be very interesting to treat such a modified scattering problem for $b = 0$. See \cite{gioz} for the Schr\"odinger case and \cite{pusa} for the semirelativistic case.
\end{rem}

We denote the \textit{charge conjugation operator} by $\mathsf C$, which is defined by
$$
\mathsf C\psi = i\gamma^2\psi^*.
$$
As usual, $\mathsf C \psi$ is interpreted as the wave function of antimatter and a direct calculation shows that $\mathsf C\psi$ satisfies the equation:
$$
-i\gamma^\mu\partial_\mu\mathsf C\psi + m\mathsf C\psi = -[V_b*(\overline{\mathsf C\psi} \gamma^5\mathsf C\psi)]\gamma^5\mathsf C\psi.
$$
The antimatter field propagates through the interaction of negative Yukawa potential.

To combine the spinor and its charge conjugation we define projection operators $P_\theta^c$ for $\theta \in \{+, -\}$ by
$$
P^c_\theta\psi = \frac12(\psi + \theta \mathsf C\psi).
$$
Then we readily get
\begin{align}\label{proj-property}
P_+^c + P_-^c = I, \quad(P^c_\theta)^2 = P^c_\theta, \quad P^c_\theta P^c_{-\theta} = 0.	
\end{align}
For this see \cite{chohlee1}. Given a spinor field $\psi:\mathbb R^{1+3}\rightarrow\mathbb C^{4}$, we further have
\begin{align}\label{charge-c-chirality}
\overline{P^c_\theta\psi} \gamma^5 P^c_\theta\psi = 0.
\end{align}
This represents the relation between charge conjugation and chirality. We append the proof in the last section.

Now $P_\theta^c\psi$ satisfy the system:
\begin{align*}
-i\gamma^\mu\partial_\mu P_+^c\psi + m P_+^c\psi = V_b*(\overline{P^c_+\psi}\gamma^5 P^c_{-}\psi+\overline{P^c_{-}\psi}\gamma^5P^c_+\psi)\gamma^5P_+^c\psi, \\
-i\gamma^\mu\partial_\mu P_-^c\psi + m P_-^c\psi = V_b*(\overline{P^c_+\psi}\gamma^5P^c_{-}\psi + \overline{P^c_{-}\psi}\gamma^5P^c_+\psi)\gamma^5P_-^c\psi.
\end{align*}
Hence, to solve the equation \eqref{main-eq} and to prove Theorem \ref{gwp}, we need to consider the system:
\begin{align}\label{c-c-decomp}
\begin{aligned}
-i\gamma^\mu\partial_\mu \varphi + m \varphi = V_b*(\overline{\varphi}\gamma^5 \chi+\overline{\chi}\gamma^5\varphi)\gamma^5\varphi, \\
-i\gamma^\mu\partial_\mu \chi + m \chi = V_b*(\overline{\varphi}\gamma^5\chi + \overline{\chi}\gamma^5\varphi)\gamma^5\chi
\end{aligned}
\end{align}
with initial data
$$
\varphi(0) = P_+^c\psi_0,\qquad \chi(0) = P^c_{-}\psi_0.
$$
If we set $\psi = P_+^c\varphi + P_-^c\chi$, then by taking $P_\pm^c$ to \eqref{c-c-decomp} one can turn back to \eqref{main-eq}. We will discuss the turning-back in the appendix.

Under the above setup we prove the following:
\begin{thm}\label{c-c-gwp}
Let $\sigma > 0$. Then there exists $0 < \epsilon < 1$ such that for any ${\mathsf A} > 0$ and any $0 < \mathtt a \le  \epsilon {\mathsf A}^{-1}$, if the initial data $\psi_0 \in L_x^{2, \sigma}$ satisfy
$$
\|P^c_+\psi_0\|_{L^{2,\sigma}} \le \mathtt a,\ \|P^c_{-}\psi_0\|_{L^{2,\sigma}}\le {\mathsf A},
$$
then the Cauchy problem of \eqref{c-c-decomp} is globally well-posed in $ C(\mathbb R;L^{2,\sigma})$.
Furthermore, there exist $\varphi^{\ell}, \chi^{\ell} \in L^{2,\sigma}$ such that $-i\gamma^\mu\partial_\mu(\varphi^\ell, \chi^\ell) + m(\varphi^\ell, \chi^\ell) = (0, 0)$
\begin{align*}
\lim_{t\rightarrow\pm\infty}\|(\varphi(t), \chi(t)) - (\varphi^\ell, \chi^\ell)\|_{L^{2,\sigma}(\mathbb R^3)} = 0.
\end{align*}
\end{thm}

To show Theorem \ref{c-c-gwp}, we reveal a null structure in $\Pi_\theta \gamma^0\gamma^5\Pi_{\theta'}$. Then we use the adapted function space consisting of $V_\theta^2$ space equipped with angular regularity, which will be introduced in Section 2. By exploiting the null structure, a contraction argument can be readily carried out by the duality as described in \cite{candyherr, candyherr1, chohlee}. We will sketch the proof in Section 3.

\section{Preliminaries}\label{sec-pre}
This section presents the preliminary setup and notations. In what follows we assume $m = 1$ for simplicity of presentation.

\subsection{Energy projection}\label{sec:dirac-op}
We further decompose the system \eqref{c-c-decomp} by energy projection $\Pi_\theta$ for $\theta\in\{+,-\}$, which is defined by
\begin{align}
\Pi_\theta(\xi) = \frac12\left(I_{4\times4}+\theta\frac{\xi_j\gamma^0\gamma^j + \gamma^0}{\Lambda(\xi)} \right),	
\end{align}
where $\Lambda(\xi) = \sqrt{1 + |\xi|^2}$. Here we used the summation convention. Now we define the Fourier multiplier by the identity $\mathcal F_x[\Pi_\theta f](\xi) = \Pi_\theta(\xi)\widehat{f}(\xi)$ and $\mathcal F_x[\Lambda(D)f](\xi) = \Lambda(\xi)\widehat f(\xi)$, where $D = -i\nabla$. By an easy computation one easily see the identity $\Pi_\theta\Pi_\theta=\Pi_\theta$ and $\Pi_\theta\Pi_{-\theta}=0$. Then we have $\psi = \sum_{\theta \in \{+, -\}}\Pi_\theta\psi$. Also we see that $\Lambda(D)(\Pi_+ - \Pi_-) = \gamma^0\gamma^j(-i\partial_j) + \gamma^0 $ and this leads us to the system: For $\theta, \theta' \in \{+, -\}$
\begin{align*}
(-i\partial_t + \theta\Lambda(D))\Pi_\theta\varphi     &= \Pi_\theta [V_b*(\overline{\varphi}\gamma^5 \chi + \overline{\chi}\gamma^5\varphi)\gamma^0\gamma^5\varphi],\\
(-i\partial_t + \theta'\Lambda(D))\Pi_{\theta'} \chi   &= \Pi_{\theta'}[V_b*(\overline{\varphi}\gamma^5 \chi + \overline{\chi}\gamma^5\varphi)\gamma^0\gamma^5\chi],\\
\varphi_\theta(0) &= \Pi_\theta[P_+^c\psi_{0}],\;\;\chi_{\theta'}(0) = \Pi_{\theta'}[P_-^c\psi_{0}].
\end{align*}
Hence solving \eqref{c-c-decomp} is equivalent to finding solutions $(\varphi_\theta, \chi_{\theta'})$ to the system of equations:
\begin{align}\begin{aligned}\label{e-decomp}
(-i\partial_t + \theta\Lambda(D))\varphi_\theta &= \sum_{\substack{\theta_j \in \{+, -\}\\j = 1, \cdots, 5}}\Pi_\theta [V_b*(\overline{\varphi_{\theta_1}}\gamma^5 \chi_{\theta_2} + \overline{\chi_{\theta_3}}\gamma^5\varphi_{\theta_4})\gamma^0\gamma^5\varphi_{\theta_5}],\\
(-i\partial_t + \theta'\Lambda(D)) \chi_{\theta'}     &=  \sum_{\substack{\theta_j' \in \{+, -\}\\j = 1, \cdots, 5}} \Pi_{\theta'}[V_b*(\overline{\varphi_{\theta_1'}}\gamma^5 \chi_{\theta_2'} + \overline{\chi_{\theta_3'}}\gamma^5\varphi_{\theta_4'})\gamma^0\gamma^5\chi_{\theta_5'}],\\
\varphi_\theta(0) &= \Pi_\theta[P_+^c\psi_{0}],\;\;\chi_{\theta'}(0) = \Pi_{\theta'}[P_-^c\psi_{0}].
\end{aligned}\end{align}
The solutions $\varphi_\theta, \chi_{\theta'}$ clearly satisfy that $\Pi_{-\theta}\varphi_\theta = \Pi_{-\theta'}\chi_{\theta'} = 0$.

\subsection{Revealing null structure}
\begin{align*}
 \gamma^0\gamma^5\Pi_\theta(\xi) &= i\gamma^1\gamma^2\gamma^3\Pi_{\theta}(\xi)\\
 & = \frac i2\gamma^1\gamma^2\gamma^3 + \frac\theta2 \Lambda(\xi)^{-1} i\gamma^1\gamma^2\gamma^3(\xi_j\gamma^0\gamma^j + \gamma^0).
\end{align*}
Using $\gamma^\mu\gamma^\nu + \gamma^\nu\gamma^\mu = 0$ for $\mu \neq \nu$ and $\gamma^j\gamma^j = -I_{4\times 4}$, we have
$$
i\gamma^1\gamma^2\gamma^3\gamma^0 = -i\gamma^0\gamma^1\gamma^2\gamma^3 = -\gamma^5 = -\gamma^0(\gamma^0\gamma^5),
$$
$$
i\gamma^1\gamma^2\gamma^3 \xi_1\gamma^0\gamma^1 = \xi_1\gamma^1\gamma^5 = -\xi_1\gamma^0\gamma^1(\gamma^0\gamma^5),
$$
$$
i\gamma^1\gamma^2\gamma^3 \xi_2\gamma^0\gamma^2 = \xi_2\gamma^2\gamma^5 = -\xi_2\gamma^0\gamma^2(\gamma^0\gamma^5),
$$
and
$$
i\gamma^1\gamma^2\gamma^3 \xi_3\gamma^0\gamma^3 = \xi_3\gamma^3\gamma^5 = -\xi_3\gamma^0\gamma^3(\gamma^0\gamma^5).
$$
Hence
$$
\gamma^0\gamma^5\Pi_\theta(\xi) = \Pi_{-\theta}\gamma^0\gamma^5.
$$
By this we deduce that
\begin{align}\label{null-gamma5}
\Pi_\theta(\xi)\gamma^0\gamma^5\Pi_{\theta'}(\eta) = \Pi_\theta(\xi)\Pi_{-\theta'}(\eta)\gamma^0\gamma^5
\end{align}
from which we expect a spinorial null structure. In fact, one can show 
\begin{align}\label{null-bound}
|\Pi_\theta(\xi)\gamma^0\gamma^5\Pi_{\theta'}(\eta)| \les {\angle}(\theta\xi, \theta'\eta) + \frac{|\theta|\xi| - \theta'|\eta| |}{\Lambda(\xi)\Lambda(\eta)},
\end{align}
where $\angle(\xi, \eta)$ is the angle between $\xi$ and $\eta$. In the scalar case the sign is positive in the second term. For this see \cite{bachelot,danfoselb, candyherr}.

\subsection{Adapted function spaces}\label{ftn-sp}
Let $\mathcal I=\left\{ \{t_k\}_{k=0}^K : t_k\in\mathbb R, t_k<t_{k+1} \right\}$ be the set of increasing sequences of real numbers.
We define the $2$-variation of $v$ to be
$$
|v|_{V^2} = \sup_{ \{t_k\}_{k=0}^K\in\mathcal I } \left( \sum_{k=0}^K\|v(t_k)-v(t_{k-1})\|_{L^2_x}^2 \right)^\frac12.
$$
Then the Banach space $V^2$ can be defined to be all right continuous functions $v:\mathbb R\rightarrow L^2_x$ such that the quantity
$$
\|v\|_{V^2} = \|v\|_{L^\infty_tL^2_x} + |v|_{V^2}
$$
is finite. Set $\|u\|_{V^2_\theta}=\|e^{\theta it\Lambda(D)}u\|_{V^2}$. We recall basic properties of $V^2_\theta$ space from \cite{candyherr, haheko, kochtavi}. In particular,  we use the following lemma to prove the scattering result.
\begin{lem}[Lemma 7.4 of \cite{candyherr}]\label{v-scatter}
	Let $u\in V^2_\theta$. Then there exists $f\in L^2_x$ such that $\|u(t)-e^{-\theta it\Lambda(D)}f\|_{L^2_x}\rightarrow0$ as $t\rightarrow\pm\infty$.
\end{lem}
\noindent We refer the readers to \cite{haheko,kochtavi} for more details.

We fix a smooth function $\rho\in C^\infty_0(\mathbb R)$ such that $\rho$ is supported in the interval $( \frac12, 2)$ and we let
$$
\sum_{\lambda\in2^{\mathbb Z}}\rho\left(\frac t\lambda\right) =1,
$$
and write $\rho_1(t) = \sum_{\lambda\le1}\rho(\frac t\lambda)$ with $\rho_1(0)=1$. Now we define the standard Littlewood-Paley multipliers, for $\lambda\in 2^{\mathbb N}$, $\lambda>1$:
$$
P_\lambda = \rho\left(\frac{|D|}{\lambda}\right),\quad P_1 = \rho_1(|D|).
$$

For a dyadic number $N>1$, we define the spherical dyadic decompositions by
\begin{align*}
H_N(f)(x) & = 	\sum_{\ell = N}^{2N-1}\sum_{m = -\ell}^{\ell}\langle f(|x|\omega), Y_{\ell m}(\omega)\rangle_{L^2_\omega(\mathbb S^2)}Y_{\ell m}\big(\frac{x}{|x|}\big), \\
H_1(f)(x) & = \sum_{\ell = 0, 1}\sum_{m = -\ell}^{\ell}\langle f(|x|\omega), Y_{\ell m}(\omega)\rangle_{L^2_\omega(\mathbb S^2)}Y_{\ell m}\big(\frac{x}{|x|}\big),
\end{align*}
where $Y_{\ell m}$ is the orthonormal spherical harmonics. Since $-\Delta_{\mathbb S^2}Y_{\ell m} = \ell(\ell+1)Y_{\ell m}$, by orthogonality one can readily get
$$\|\Lambda_{\mathbb S^2}^\sigma f\|_{L^2_\omega({\mathbb S^2})} \approx \left\|\sum_{N \in 2^{\mathbb N\cup\{0\}}}N^\sigma H_Nf\right\|_{L^2_\omega({\mathbb S^2})}.$$


\section{Trilinear estimates}
We define the Banach space associated with the adapted space to be the set
$$
F_\theta^{\sigma} = \left\{ \phi \in C(\mathbb R; L_x^{2, \sigma}): \Pi_{-\theta}\phi = 0\;\;\mbox{and}\;\; \|\phi\|_{ F_\theta^{\sigma} }< \infty \right\} \;\;\mbox{for }\;\;\theta \in \{+, -\},
$$
where the norm is defined by
\begin{align*}
\|\phi\|_{F_{\theta}^\sigma} :=  \bigg( \sum_{{\rm dyadic}\;\lambda, N \ge 1}N^{2\sigma}\|  P_\lambda H_N  \phi\|_{V^2_{\theta}}^2 \bigg)^\frac12.
\end{align*}
Note that $\sum_{N \ge 1}\Pi_{-\theta}H_N\phi = 0$ for all $\phi \in F_\theta^\sigma$ (see \cite{chohlee}).

A crucial part of the proof of Theorem \ref{c-c-gwp} is the following trilinear estimate: Let $\sigma>0$ and $\theta, \theta_j \in \{+, -\}$. Then for any $\varphi \in F_{\theta_1}^\sigma$, $\chi \in F_{\theta_2}^\sigma$, and $\psi \in F_{\theta_3}^\sigma$ there holds
\begin{align}\begin{aligned}\label{tri-main}
\|\mathfrak J^\theta[\Pi_\theta (V_b*(\varphi^\dagger \gamma^0\gamma^5 \chi)\gamma^0\gamma^5\psi)]\|_{ F^{\sigma}_\theta} & \lesssim \|\varphi\|_{ F_{\theta_1}^{\sigma}}	\|\chi\|_{ F_{\theta_2}^{\sigma}}\|\psi\|_{ F_{\theta_5}^{\sigma}}.
\end{aligned}\end{align}
 Here $\mathfrak J^\theta[F]$ is the Duhamel integral which reads
$$
\mathfrak J^\theta[F] = \int_0^t e^{-\theta i(t-t')\Lambda(D)}F(t')\,dt'.
$$

Now we introduce the following frequency-localised bilinear estimates.
\begin{lem}\label{main-bi}
Let $\epsilon>0$. Then there exists $\delta > 0$ such that the following frequency-localised $L^2$-bilinear estimates hold:
\begin{align*}
	&\|P_{\mu}H_N((\Pi_{\theta_1}P_{\lambda_1}H_{N_1}\varphi)^\dagger \gamma^0\gamma^5\Pi_{\theta_2}P_{\lambda_2}H_{N_2}\chi\|_{L^2_tL^2_x} \\
	 & \qquad\qquad \lesssim \mu\left(\frac{\min\{\mu,\lambda_1,\lambda_2\}}{\max\{\mu,\lambda_1,\lambda_2\}} \right)^\mathtt \delta (\min\{N_1,N_2\})^\epsilon \|P_{\lambda_1}H_{N_1}\varphi\|_{V^2_{\theta_1}}\|P_{\lambda_2}H_{N_2}\chi\|_{V^2_{\theta_2}}
\end{align*}	
for any dyadic numbers $\mu,\lambda_1,\lambda_2, N,N_1,N_2\ge1$.
\end{lem}
In view of the null condition \eqref{null-gamma5} and null bound \eqref{null-bound}, the bilinear form of Lemma \ref{main-bi} has the same type null structure as the bilinear form without $\gamma^5$. Hence one can carry out a similar way of proof for Lemma \ref{main-bi}  to the one of Proposition 3.1 of \cite{chohlee} by using Lemma \ref{main-bi} together with modulation estimates and angular estimates \cite{candyherr, chohlee, cholee}.  Indeed, we have for $\mu\lesssim\lambda_1\approx\lambda_2$,
 $$
 \|P_\mu(\varphi_{\lambda_1}^\dagger\gamma^0\gamma^5\psi_{\lambda_2})\|_{L^2_{t,x}}\lesssim\mu\|\varphi_{\lambda_1}\|_{V^2_{\theta_1}}\|\psi_{\lambda_2}\|_{V^2_{\theta_2}},
 $$
 especially when $\theta_1=\theta_2$.
 For the proof, see Proposition 3.7 of \cite{cyang}. On the other hand, we exploit an extra weight of angular regularity to get the bound as
 for some $\mathfrak d>0$,
\begin{align*}
&\|P_\mu H_N(\varphi_{\lambda_1,N_1}^\dagger\gamma^0\gamma^5\psi_{\lambda_2,N_2})\|_{L^2_tL^2_x} \\
& \lesssim \mu\left(\frac\mu{\min\{\lambda_1,\lambda_2\}}\right)^{\mathfrak d}	\min\{N_1,N_2\}\|\varphi_{\lambda_1,N_1}\|_{V^2_{\theta_1}}\|\psi_{\lambda_2,N_2}\|_{V^2_{\theta_2}}.
\end{align*}
For this we refer to \cite{hong}.
   Hence we simply combine two bound to get
   \begin{align*}
&\|P_\mu H_N(\varphi_{\lambda_1,N_1}^\dagger\gamma^0\gamma^5\psi_{\lambda_2,N_2})\|_{L^2_tL^2_x} \\
& \lesssim \mu\left(\frac\mu{\min\{\lambda_1,\lambda_2\}}\right)^{\frac\delta{8}}	(\min\{N_1,N_2\})^\delta\|\varphi_{\lambda_1,N_1}\|_{V^2_{\theta_1}}\|\psi_{\lambda_2,N_2}\|_{V^2_{\theta_2}},	
\end{align*}
for an arbitrarily small $\delta\ll1$, where we write $\psi_{\lambda,N}:=P_\lambda H_N \psi$, for brevity.

By Lemma \ref{main-bi} and the duality argument for the proof of Proposition 3.1 of \cite{chohlee}, we obtain the desired trilinear estimate \eqref{tri-main}.


\section{Proof of Theorem \ref{c-c-gwp}}

For $\theta, \theta' \in \{+, -\}$ we consider the set
 \begin{align*}
 \pmb{ \mathfrak F} = \{ \pmb \phi := &(\varphi_{+}, \varphi_{-}, \chi_{+}, \chi_- ) \in \pmb F^\sigma := F_{+}^{\sigma} \times F_{-}^{\sigma} \times F_{+}^{\sigma} \times F_{-}^{\sigma} : \\
 &\|\varphi_{\theta}\|_{F_{\theta}^{\sigma}} \le 2\|\varphi_{\theta}(0)\|_{L^{2,\sigma}},\ \|\chi_{\theta'}\|_{F_{\theta'}^{\sigma}} \le 2\|\chi_{\theta'}(0)\|_{L^{2,\sigma}} \}
 \end{align*}
 and for ${\mathsf A}, \mathtt a > 0$, we define the norm
 $$
 \|\pmb \phi \|_{\pmb{ \mathfrak F}} := \mathtt a^{-1}\sum_{\theta \in \{+, -\}}\|\varphi_{\theta}\|_{F_{\theta}^{\sigma}} + {\mathsf A}^{-1}\sum_{\theta'\in \{+, -\}}\|\chi_{\theta'}\|_{F_{\theta'}^{\sigma}}.
 $$
 Then $\pmb{\mathfrak{F}}$ is a complete metric space with the metric derived by the norm. Now we let $\pmb{\mathfrak M} = (\Phi_+, \Phi_-, X_+, X_-)$ be the inhomogeneous solution map on $\pmb{ \mathfrak F}$ for \eqref{e-decomp} given by the Duhamel's principle. Then the definition of the set $\pmb{ \mathfrak F}$ and trilinear estimate \eqref{tri-main} give us the estimate: For any $\pmb \phi \in \pmb{ \mathfrak F}$,
 \begin{align}
 \begin{aligned}
 &\sum_{\theta}\|\Phi_\theta(\pmb\phi)\|_{F_{\theta}^{\sigma}}\\
  & \le \sum_\theta\|\varphi_{\theta}(0)\|_{L^{2,\sigma}} + C (\sum_{\theta_1 \in \{+, -\}}\|\varphi_{\theta_1}\|_{F_{\theta_1}^{\sigma}})^2(\sum_{\theta_2 \in \{+, -\}}\|\chi_{\theta_2}\|_{F^{\sigma}}) \\
 & \le \sum_\theta\|\varphi_\theta(0)\|_{L^{2,\sigma}}	+16C(\sum_{\theta_1} \|\varphi_{\theta_1}(0)\|_{L^{2,\sigma}})^2(\sum_{\theta_2}\|\chi_{\theta_2}(0)\|_{L^{2,\sigma}}) \\
 & \le 	(1+16C{\mathsf A}\mathtt a)\sum_{\theta}\|\varphi_{\theta}(0)\|_{L^{2,\sigma}},
 \end{aligned}
 \end{align}
and also
\begin{align}
\begin{aligned}
&\sum_{\theta'}\|X_{\theta'}(\pmb\phi)\|_{F_{\theta'}^{\sigma}}\\
 & \le \sum_\theta\|\chi_{\theta'}(0)\|_{L^{2,\sigma}} + C (\sum_{\theta_1 \in \{+, -\}}\|\varphi_{\theta_1}\|_{F_{\theta_1}^{\sigma}})(\sum_{\theta_2 \in \{+, -\}}\|\chi_{\theta_2}\|_{F^{\sigma}})^2 \\
 & \le \sum_\theta\|\varphi_\theta(0)\|_{L^{2,\sigma}}	+ 16C(\sum_{\theta_1} \|\varphi_{\theta_1}(0)\|_{L^{2,\sigma}})(\sum_{\theta_2}\|\chi_{\theta_2}(0)\|_{L^{2,\sigma}})^2 \\
 & \le 	(1 + 16C{\mathsf A}\mathtt a)\sum_{\theta'}\|\chi_{\theta'}(0)\|_{L^{2,\sigma}}.
\end{aligned}	
\end{align}
Then we put $\mathtt a \le \dfrac{1}{16C{\mathsf A}}$ and deduce that the map $\pmb{\mathfrak M}$ is the flow map from $\pmb{ \mathfrak F}$ onto $\pmb{ \mathfrak F}$. The trilinear estimate \eqref{tri-main} leads us that the map $\pmb{\mathfrak M}$ is a contraction on the set $\pmb{ \mathfrak F}$. Indeed, suppose that we have $\pmb \phi^1, \,\pmb\phi^2 \in \pmb{ \mathfrak F}$. Then we estimate
\begin{align*}
&\sum_{\theta}\|\Phi_\theta(\pmb \phi^1) - \Phi_\theta(\pmb \phi^2)\|_{F_\theta^\sigma} \\
&\le 16C{\mathsf A}\mathtt a \sum_{\theta}\| \phi_\theta^1 -  \phi_\theta^2\|_{F_\theta^\sigma} + 8C\mathtt a^2\sum_{\theta'}\|\chi_{\theta'}^1 -\chi_{\theta'}^2\|_{F_{\theta'}^{\sigma}}
\end{align*}
and
\begin{align*}
&\sum_{\theta'}\|X_{\theta'}(\pmb \phi^1) - X_{\theta'}(\pmb \phi^2)\|_{F_{\theta'}^{\sigma}}\\
& \le 16C{\mathsf A}\mathtt a \sum_{\theta'}\|\chi_{\theta'}^1 - \chi_{\theta'}^2\|_{F_{\theta'}^{\sigma}} + 8C{\mathsf A}^{2}\sum_{\theta'}\|\varphi_{\theta'}^1 - \varphi_{\theta'}^2\|_{F_{\theta'}^{\sigma}}.
\end{align*}
In consequence we obtain
\begin{align*}
\|\pmb{\mathfrak M}(\pmb \phi^1) - \pmb{\mathfrak M}(\pmb \phi^2)\|_{\pmb{ \mathfrak F}} &\le 24C{\mathsf A}\sum_{\theta}\|\varphi_{\theta}^1 -\varphi_{\theta}^2\|_{F_{\theta}^{\sigma}} + 24C\mathtt a\sum_{\theta'}\|\chi_{\theta'}^1-\chi_{\theta'}^2\|_{F_{\theta'}^{\sigma}}\\
& = 24C{\mathsf A}\mathtt a\|\pmb \phi^1 - \pmb\phi^2\|_{\pmb{ \mathfrak F}}\,\,.
\end{align*}
Thus by choosing $\epsilon = \dfrac{1}{48C}$, the solution map $\pmb{\mathfrak M}$ is a contraction on $\pmb{ \mathfrak F}$ for any $\mathtt a \le \epsilon{\mathsf A}^{-1}$. The scattering follows immediately from Lemma \ref{v-scatter}.

\section{Appendix}
\subsection{Proof of \eqref{charge-c-chirality}}
We prove \eqref{charge-c-chirality}. By a direct calculation one can readily show that $\overline{P^c_\theta\psi}\gamma^5P^c_\theta\psi$ is purely imaginary. We now write
\begin{align}\label{sum-c-c-c}
&\overline{P^c_\theta\psi}\gamma^5 P^c_\theta\psi \nonumber\\
 &= \frac14(\psi^\dagger + \theta \psi^T(i\gamma^2))\gamma^0\gamma^5(\psi+\theta i\gamma^2\psi^*),\tag{A1}\\
&= \psi^\dagger \gamma^0\gamma^5 \psi + \psi^\dagger \gamma^0\gamma^5 (i\gamma^2)\psi^* + \psi^T(i\gamma^2)\gamma^0\gamma^5\psi + \psi^T(i\gamma^2)\gamma^0\gamma^5(i\gamma^2)\psi^*.\nonumber	
\end{align}
Then by the relation $\gamma^\mu\gamma^\nu + \gamma^\nu\gamma^\mu = 0 (\mu \neq \nu)$ and $\gamma^2\gamma^2 = -I_{4\times 4}$ we have
$$
\psi^T(i\gamma^2)\gamma^0\gamma^5(i\gamma^5)\psi^* = \psi^T\gamma^0\gamma^5\psi^*.
$$
Since $(\psi^\dagger\gamma^0\gamma^5\psi)^* = \psi^T\gamma^0\gamma^5\psi^*$, the sum of the first and last terms in \eqref{sum-c-c-c} are real-valued.
By definition of $\gamma^2$, $(\gamma^2)^* = -\gamma^2$ and hence the sum of the second and third terms in \eqref{sum-c-c-c} are also real-valued.
The LHS of \eqref{sum-c-c-c} is purely imaginary, whereas the RHS is purely real. Therefore $\overline{P^c_\theta\psi}\gamma^5 P^c_\theta\psi = 0$.

\subsection{From \eqref{c-c-decomp} to \eqref{main-eq}}

Let $(\varphi, \chi)$ be the solution of \eqref{c-c-decomp} satisfying the condition of Theorem \ref{c-c-gwp}.
Taking $P_-^c$ to the first equation of \eqref{c-c-decomp}, we have the equation
$$
-i\gamma^\mu \partial_\mu P_-^c\varphi + m P_-^c\varphi = V_b*(\overline{\varphi}\gamma^5 \chi+\overline{\chi}\gamma^5\varphi)\gamma^5P_-^c\varphi.
$$
Hence $P_-^c\varphi$ is the solution with initial data $P_-^c\varphi(0) = P_-^cP_+^c\psi(0) = 0$ and it can be written as
$$
\Pi_\theta (P_-^c\varphi) = i\int_0^t e^{-\theta i(t-t')\Lambda(D)}\Pi_\theta[V_b*(\overline{\varphi}\gamma^5\chi + \overline{\chi}\gamma^5\varphi)\gamma^0\gamma^5 P_-^c\varphi]\,dt'.
$$
By trilinear estimates \eqref{tri-main} and the choice of $\mathtt a$ we have
$$
\sum_{\theta}\|\Pi_\theta P_-^c\varphi\|_{F_\theta^\sigma} \le 16C\mathtt a\mathsf A \sum_{\theta}\|\Pi_\theta P_-^c\varphi\|_{F_\theta^\sigma} \le \frac12\sum_{\theta}\|\Pi_\theta P_-^c\varphi\|_{F_\theta^\sigma}.
$$
Therefore $P_-^c\varphi = 0$. In the same way we deduce that $P_+^c\chi = 0$. We now write the system \eqref{c-c-decomp} as
\begin{align}\label{c-c-decomp-1}
-i\gamma^\mu\partial_\mu \varphi + m \varphi = V_b*(\overline{P_+^c\varphi}\gamma^5 P_-^c\chi+\overline{P_-^c\chi}\gamma^5P_+^c\varphi)\gamma^5\varphi,\nonumber \\
-i\gamma^\mu\partial_\mu \chi + m \chi = V_b*(\overline{P_+^c\varphi}\gamma^5P_-^c\chi + \overline{P_-^c\chi}\gamma^5P_+^c\varphi)\gamma^5\chi.\tag{A2}
\end{align}
Taking $P_+^c$ and $P_-^c$ to \eqref{c-c-decomp-1}, we finally get
\begin{align*}
-i\gamma^\mu\partial_\mu P_+^c\varphi + m P_+^c\varphi = V_b*(\overline{P_+^c\varphi}\gamma^5 P_-^c\chi+\overline{P_-^c\chi}\gamma^5P_+^c\varphi)\gamma^5P_+c\varphi, \\
-i\gamma^\mu\partial_\mu P_-^c\chi + m P_-^c\chi = V_b*(\overline{P_+^c\varphi}\gamma^5P_-^c\chi + \overline{P_-^c\chi}\gamma^5P_+^c\varphi)\gamma^5P_-c\chi.
\end{align*}
By setting $\psi = P_+\varphi + P_-\chi$ we obtain the original equation \eqref{main-eq}.
\section*{Acknowledgements}
Y. Cho was supported by NRF-2021R1I1A3A04035040(Republic of Korea).


\section*{Data Availability}
The data that support the findings of this study are available within the article.



\end{document}